\documentclass[twoside,11pt,reqno]{amsart}

\usepackage{amsmath,amsthm,amssymb,amstext,amsfonts,amscd}
\usepackage{graphicx}
\usepackage{multirow}
\usepackage{url}
 \newtheorem{theorem}{Theorem}[section]

\numberwithin{equation}{section}

\begin{document}

\title[Certain study of Generalized Apostol-Bernoulli poly-Daehee poly...]
{Certain study of Generalized Apostol-Bernoulli poly-Daehee polynomials and its properties}

\author[{\bf N. U. Khan and S. Husain}]{\bf Nabiullah Khan and Saddam Husain}

\address{Nabiullah Khan: Department of Applied
Mathematics, Faculty of Engineering and Technology,
      Aligarh Muslim University, Aligarh 202002, India}
 \email{nukhanmath@gmail.com}

\bigskip
\address{Saddam Husain: Department of Applied
    Mathematics, Faculty of Engineering and Technology,
    Aligarh Muslim University, Aligarh 202002, India}
 \email{saddamhusainamu26@gmail.com}



\begin{abstract}
In this paper, we present a new type of generating function of generalized Apostol-Bernoulli poly Daehee polynomial (GABPDP). By using generating function of GABPDP, we discus some special cases and useful identites of generalized Apostol-Bernoulli poly Daehee polynomials. We also drive implicit summation formulae of it.
\end{abstract}

\maketitle

\section{\bf{Introduction}}
In the present article, we take following useful standard notation: $\mathbb{N}_0$=: $\mathbb{N}\cup {(0)}$, $\mathbb{Z}^{-}$, $\mathbb{N}$, $\mathbb{R}$ and $\mathbb{C}$ be set of all negative integer, natural, real and complex number.\\
We know well known classical poly-logarithm function $L_{ik}(z)$ (see \cite{khan-usman-aman,wasim}) is given by

\begin{equation}\label{e1.1}
L_{ik}(z)= \sum_{n=0}^{\infty} \frac{z^m}{m^k},\,\,\,\,\,\,\,\,\,\,(z\in \mathbb{C}, k\in \mathbb{Z}).
\end{equation}

For $k=1$, $L_{ik}(z)=-log(1-z)$.\\

Recently, many authors (see \cite{carlitz,choi-nu-usman1,jolany-mohsen,khan-usman-aman,khan-usman1,khan-usman2,khan-usman3,khan-usman-choi2,subuhi-tabinda,kim-taekyun,dae-taekyun,Takao-luca})  have studied the Daehee polynomial, Bernoulli polynomial and Euler polynomial. Also studied second kind Bernoulli polynomial, poly Bernoulli polynomial and generalized Bernoulli, Euler and Genocchi polynomials. Which is defined as follows:\\

The $n^{th}$ Daehee polnomial and Daehee number (see \cite{lim-kwon}) are defined as follows:
\begin{equation}\label{e1.2}
\sum_{n=0}^{\infty}\frac{x^n}{n!}{\mathfrak{D}}_{n}(\gamma)=(1+x)^{\gamma}\,\,\frac{log(1+x)}{x}.
\end{equation}

If we take $\gamma = 0$ then ${\mathfrak{D}}_{n}:= {\mathfrak{D}}_{n}(0)$ are the Daehee number.\\

The Euler polynomial ${\mathcal{E}}_{n}(\gamma)$ and Euler number ${\mathcal{E}}_{n}(0)$ are defined by the following generating function to be

\begin{equation}\label{e1.3}
e^{\gamma{x}}\left(\frac{2}{e^{x}+1}\right)=\sum_{n=0}^{\infty}\frac{x^n}{n!}~~{\mathcal{E}}_{n}(\gamma), \,\,\,\,\,\,\, ( |x|<2\pi ).
\end{equation}

If we take $\gamma = 0$ then the Euler number $\mathcal{E}_{n}(0):= \mathcal{E}_{n}$ are defined by

\begin{equation}\label{e1.4}
\left(\frac{2}{e^{x}+1}\right)=\sum_{n=0}^{\infty}\frac{x^n}{n!}~\mathcal{E}_{n}, \,\,\,\,\,\,\, ( |x|<2\pi ).
\end{equation}

The Bernoulli polynomial and Bernoulli number is introduced by Datolli et al. (see\cite{Dattoli}), which is defined by the following generating function to be

\begin{equation}\label{e1.5}
e^{\gamma{x}}\left(\frac{x}{e^{x}-1}\right)=\sum_{n=0}^{\infty}\frac{x^n}{n!}{\mathfrak{B}}_{n}(\gamma), \,\,\,\,\,\,\, ( |x|<2\pi ).
\end{equation}

If we take $\gamma = 0$ then the Euler number ${\mathfrak{B}}_{n}(0):= {\mathfrak{B}}_{n}$ are defined by

\begin{equation}\label{e1.6}
\left(\frac{x}{e^{x}-1}\right)=\sum_{n=0}^{\infty}\frac{x^n}{n!}{\mathfrak{B}}_{n}, \,\,\,\,\,\,\, ( |x|<2\pi ).
\end{equation}

The poly-Bernoulli polynomial ${\mathfrak{B}}_{n}^{(k)}(\gamma)$ are defined by the following generating function to be

\begin{equation}\label{e1.7}
e^{\gamma{x}}\left(\frac{L_{ik}(1-e^{-x})}{e^{x}-1}\right)=\sum_{n=0}^{\infty}\frac{x^n}{n!}~{\mathfrak{B}}_{n}^{(k)}(\gamma).
\end{equation}

If we take $\gamma=0$, then poly-Bernoulli number ${\mathfrak{B}}_{n}^{(k)}(0):={\mathfrak{B}}_{n}^{(k)}$ are defined by

\begin{equation}\label{e1.8}
\frac{L_{ik}(1-e^{-x})}{e^{x}-1}=\sum_{n=0}^{\infty}\frac{x^n}{n!}~{\mathfrak{B}}_{n}^{(k)}.
\end{equation}

\vspace{0.15cm}
For $k=1$ in \eqref{e1.7} and \eqref{e1.8}, we get ${\mathfrak{B}}_{n}^{(1)}(\gamma):={\mathfrak{B}}_{n}(\gamma)$,  ${\mathfrak{B}}_{n}^{(1)}:={\mathfrak{B}}_{n}$.

The higher order poly-Bernoulli polynomial ${\mathfrak{B}}_{n}^{(k, r)}(\gamma)$ are defined by the following generating function to be

\begin{equation}\label{e1.7}
e^{\gamma{x}}\left(\frac{L_{ik}(1-e^{-x})}{e^{x}-1}\right)^{r}=\sum_{n=0}^{\infty}\frac{x^n}{n!}~{\mathfrak{B}}_{n}^{(k, r)}(\gamma).
\end{equation}

The poly-Daehee polynomials ${\mathfrak{D}}_{n}^{(k)}(\gamma)$ (see \cite{lim-kwon}) are defined by 
\begin{equation}\label{e1.9}
\frac{\log(1+x)}{L_{ik}(1-e^{-x})}(1+x)^{\gamma}=\sum_{n=0}^{\infty}{\mathfrak{D}}_{n}^{(k)}(\gamma)\frac{x^{n}}{n!}\quad(k\in \mathbb{Z_{+}}=\mathbb{N}\cup \{0\}).
\end{equation}
If $\gamma=0, {\mathfrak{D}}_{n}^{(k)}={\mathfrak{D}}_{n}^{(k)}(0)$ are called the poly-Daehee numbers, which are defined as
\begin{equation}\label{e1.10}
\frac{\log(1+x)}{L_{ik}(1-e^{-x})}=\sum_{n=0}^{\infty}{\mathfrak{D}}_{n}^{(k)}\frac{x^{n}}{n!}\quad(k\in \mathbb{Z_{+}}=\mathbb{N}\cup \{0\}).
\end{equation}

\vspace{0.15cm}
The poly-Bernoulli polynomial $b_{n}^{(k)}(\gamma)$ of second kind (see \cite{kim-kwon-lee}) are defined by the following generating function

\begin{equation}\label{e1.11}
(1+x)^{\gamma}\left(\frac{L_{ik}(1-e^{-x})}{log(1+x)}\right)=\sum_{n=0}^{\infty}\frac{x^n}{n!}~\mathfrak{b}_{n}^{(k)}(\gamma).
\end{equation}

If  $\gamma=0$, then  $\mathfrak{b}_{n}^{(k)}(0):=\mathfrak{b}_{n}^{(k)}$ are called poly-Bernoulli number of second kind.\\

The generalized Bernoulli, Euler and Genocchi polynomials of order $a\in \mathbb{C}$ and $|x|< 2\pi$ are defined by the following generating relation respectively;

\begin{equation}\label{e1.12}
e^{\gamma{x}}\left(\frac{x}{e^{x}-1}\right)^{a}=\sum_{n=0}^{\infty}\frac{x^n}{n!}{\mathfrak{B}}_{n}^{(a)}(\gamma), 
\end{equation}

\vspace{0.25cm}
\begin{equation}\label{e1.13}
e^{\gamma{x}}\left(\frac{2}{e^{x}+1}\right)^{a}=\sum_{n=0}^{\infty}\frac{x^n}{n!}\mathcal{E}_{n}^{(a)}(\gamma), 
\end{equation}

\vspace{0.25cm}
\begin{equation}\label{e1.12}
e^{\gamma{x}}\left(\frac{2x}{e^{x}+1}\right)^{a}=\sum_{n=0}^{\infty}\frac{x^n}{n!}~{\mathcal{G}}_{n}^{(a)}(\gamma). 
\end{equation}

\vspace{0.15cm}
If $\gamma = 0$ in generating function \eqref{e1.10}, \eqref{e1.11} and \eqref{e1.12}, we get

\begin{equation}\label{e1.14}
\left(\frac{x}{e^{x}-1}\right)^{a}=\sum_{n=0}^{\infty}\frac{x^n}{n!}~{\mathfrak{B}}_{n}^{(a)}(0), 
\end{equation}

\vspace{0.25cm}
\begin{equation}\label{e1.15}
\left(\frac{2}{e^{x}+1}\right)^{a}=\sum_{n=0}^{\infty}\frac{x^n}{n!}~\mathcal{E}_{n}^{(a)}(0), 
\end{equation}

\vspace{0.25cm}
\begin{equation}\label{e1.16}
\left(\frac{2x}{e^{x}+1}\right)^{a}=\sum_{n=0}^{\infty}\frac{x^n}{n!}~\mathcal{G}_{n}^{(a)}(0),
\end{equation}
where
\begin{equation*}
{\mathfrak{B}}_{n}^{(a)}(0):={\mathfrak{B}}_{n}^{(a)},~~ \mathcal{E}_{n}^{(a)}(0):=\mathcal{E}_{n}^{(a)},~~ \mathcal{G}_{n}^{(a)}(0):=\mathcal{G}_{n}^{(a)}\,\,\,\,\,(n\in \mathbb{N}_0),
\end{equation*}

\vspace{0.15cm}
are the generalized Bernoulli, Euler and Genocchi numbers respectively.\\

It can be seen that if $a=1$ in \eqref{e1.10}, \eqref{e1.11} and \eqref{e1.12} respectively, then 
\begin{equation*}
{\mathfrak{B}}_{n}^{(1)}(\gamma):={\mathfrak{B}}_{n}(\gamma),~~  \mathcal{E}_{n}^{(1)}(\gamma):=\mathcal{E}_{n}(\gamma),~~  \mathcal{G}_{n}^{(1)}(\gamma):=\mathcal{G}_{n}(\gamma)\,\,\,\,\,(n\in \mathbb{N}_0).
\end{equation*}

\vspace{0.15cm}
The generalized Apostol-Bernoulli polynomial $\mathfrak{B}_{n}^{a}(\gamma, \lambda)$ of order $a \in \mathbb{C}$ (see \cite{Luo-ming-sri}), defined by the following generating function to be

\begin{equation}\label{e1.17}
e^{\gamma{x}}\left(\frac{x}{\lambda e^{x}-1}\right)^{a}=\sum_{n=0}^{\infty}\frac{x^n}{n!}~{\mathfrak{B}}_{n}^{(a)}(\gamma; \lambda),\,\,\,\,(|x+ln{\lambda}|<2\pi),
\end{equation}
with
$${\mathfrak{B}}_{n}^{(a)}(\gamma; 1):={\mathfrak{B}}_{n}^{(a)}(\gamma)$$
and
$${\mathfrak{B}}_{n}^{(a)}(0; \lambda):={\mathfrak{B}}_{n}^{(a)}(\lambda)$$

which is known as Apostol-Bernoulli number ${\mathfrak{B}}_{n}^{(a)}(\lambda)$ of order $a$.

\vspace{0.15cm}
The generalized Apostol-Euler polynomial $\mathcal{E}_{n}^{a}(\gamma, \lambda)$ of order $a \in \mathbb{C}$ (see\cite{luo-ming}), defined by the following generating function to be

\begin{equation}\label{e1.18}
e^{\gamma{x}}\left(\frac{2}{\lambda e^{x}+1}\right)^{a}=\sum_{n=0}^{\infty}\frac{x^n}{n!}~\mathcal{E}_{n}^{(a)}(\gamma; \lambda),\,\,\,\,(|x+ln{\lambda}|<2\pi),
\end{equation}
with
$$\mathcal{E}_{n}^{(a)}(\gamma; 1):=\mathcal{E}_{n}^{(a)}(\gamma)$$
and
$$\mathcal{E}_{n}^{(a)}(0; \lambda):=\mathcal{E}_{n}^{(a)}(\lambda)$$

which is known as Apostol-Euler number of order a.

\vspace{0.15cm}
The generalized Apostol-Genocchi polynomial $\mathcal{G}_{n}^{a}(\gamma, \lambda)$  of order $a \in \mathbb{C}$ (see \cite{luo-ming2}), defined by the following generating function to be

\begin{equation}\label{e1.19}
e^{\gamma{x}}\left(\frac{2x}{\lambda e^{x}+1}\right)^{a}=\sum_{n=0}^{\infty}\frac{x^n}{n!}~\mathcal{G}_{n}^{(a)}(\gamma; \lambda),\,\,\,\,(|x+ln{\lambda}|<2\pi),
\end{equation}
with
$$\mathcal{G}_{n}^{(a)}(\gamma; 1):=\mathcal{G}_{n}^{(a)}(\gamma)$$
and
$$\mathcal{G}_{n}^{(a)}(0; \lambda):=\mathcal{G}_{n}^{(a)}(\lambda)$$

which is known as Apostol-Genocchi number of order a.\\

The generalized Apostol-Bernoulli polynomial $\mathfrak{B}_{n, a}^{[m-1]}(\eta, \lambda)$ of order $a \in \mathbb{C}$ (see \cite{Tremblay-Gaboury}), defined by the following generating function to be
\begin{equation}\label{e1.22}
e^{\eta{x}}\left(\frac{x^m}{\lambda e^{x}-\sum\limits_{l=0}^{m-1}\frac{x^l}{l!}}\right)^{a}=\sum_{n=0}^{\infty}\frac{x^n}{n!}~{\mathfrak{B}}_{n, a}^{[m-1]}(\eta; \lambda)
\end{equation}

If $\eta=0$, the generalized Apostol-Bernoulli number defined by:
\begin{equation}\label{e1.23}
\left(\frac{x^m}{\lambda e^{x}-\sum\limits_{l=0}^{m-1}\frac{x^l}{l!}}\right)^{a}=\sum_{n=0}^{\infty}\frac{x^n}{n!}~{\mathfrak{B}}_{n, a}^{[m-1]}(\lambda)
\end{equation}

\section{\bf{Generalized Apostol-Bernoulli Poly-Daehee Polynomials}}

In this section, we define a generalization and unification of generalized Apostol-Bernoulli poly-Daehee polynomials (GABPDP) and define a certain useful properties and implicit formulae of Apostol-Bernoulli poly-Daehee Polynomials. which is defined as follows:

\vspace{0.10cm}
For $\gamma, \eta \in \mathbb{R}, n, m \in \mathbb{N}$ and $a \in \mathbb{C}$ we have define the following generating function for generalized Apostol-Bernoulli poly Daehee polynomials (GABPDP):

\begin{equation}\label{e2.1}
(1+x)^{\gamma}\,\,\frac{log(1+x)}{L_{ik}(1-e^{-x})} e^{\eta{x}}\left(\frac{x^m}{\lambda e^{x}-\sum\limits_{l=0}^{m-1}\frac{x^l}{l!}}\right)^{a}=\sum_{n=0}^{\infty}\frac{x^n}{n!}\, {_{\mathfrak{B}} {\mathfrak{D}}_{n, a}^{[k, m-1]}(\gamma, \eta; \lambda)}.
\end{equation}

If $\gamma = \eta = 0$, then $_{\mathfrak{B}} {\mathfrak{D}}_{n, a}^{[k, m-1]}(0, 0; \lambda) :=  {_{\mathfrak{B}} {\mathfrak{D}}_{n, a}^{[k, m-1]}(\lambda)}$ are called generalized Apostol-Bernoulli poly-Daehee number (GABPDN).

\subsection{\bf Special Cases}
In this section, we discus some particular cases of generalized Apostol-Bernoulli poly-Daehee polynomials:
\begin{enumerate}
	\item If $m=1$, then equation \eqref{e2.1} become Apostol-Bernoulli poly-Daehee polynomials with generating function
	\begin{equation*}
	(1+x)^{\gamma}\,\,\frac{log(1+x)}{L_{ik}(1-e^{-x})} e^{\eta{x}}\left(\frac{x}{\lambda e^{x}-1}\right)^{a}=\sum_{n=0}^{\infty}\frac{x^n}{n!}\, {_{\mathfrak{B}} {\mathfrak{D}}_{n, a}^{(k)}(\gamma, \eta; \lambda)}.
	\end{equation*}
	\item If $m=1, \lambda=1$ and $a=1$, then equtaion \eqref{e2.1} become Bernoulli poly-Daehee polynomials with generating function
	\begin{equation*}
	(1+x)^{\gamma}\,\,\frac{log(1+x)}{L_{ik}(1-e^{-x})} e^{\eta{x}}\left(\frac{x}{ e^{x}-1}\right)=\sum_{n=0}^{\infty}\frac{x^n}{n!}\, {_{\mathfrak{B}} {\mathfrak{D}}_{n}^{(k)}(\gamma, \eta)}.
	\end{equation*}
	\item If $m=1, k=1$, then equation \eqref{e2.1} become Bernoulli based Daehee polynomials with generating function
	\begin{equation*}
	(1+x)^{\gamma}\,\,\frac{log(1+x)}{x} e^{\eta{x}}\left(\frac{x}{\lambda e^{x}-1}\right)^{a}=\sum_{n=0}^{\infty}\frac{x^n}{n!}\, {_{\mathfrak{B}} {\mathfrak{D}}_{n, a}(\gamma, \eta; \lambda)}.
	\end{equation*}
	\item If $m=1, \eta=0$ and $a=0$, then \eqref{e2.1} become poly Daehee polynomials with generating function
	\begin{equation*}
	(1+x)^{\gamma}\,\,\frac{log(1+x)}{L_{ik}(1-e^{-x})} =\sum_{n=0}^{\infty}\frac{x^n}{n!}\, { {\mathfrak{D}}_{n}^{(k)}(\gamma)}.
	\end{equation*}
	\item If $m=1, k=1, \eta=0$ and $a=0$, then \eqref{e2.1} become $n^{th}$ Daehee polynomials with generating function
	\begin{equation*}
	(1+x)^{\gamma}\,\,\frac{log(1+x)}{x} =\sum_{n=0}^{\infty}\frac{x^n}{n!}\, { {\mathfrak{D}}_{n}(\gamma)}.
	\end{equation*}
\end{enumerate}

\begin{theorem}
Let $\gamma, \eta \in \mathbb{R}$, $n, m \in \mathbb{N}_0$ then GABPDP satisfy the following relation:
\begin{equation}\label{e2.2}
{_{\mathfrak{B}} {\mathfrak{D}}_{n, a}^{[k, m-1]}(\gamma, \eta; \lambda)} = \sum\limits_{j=0}^{n}\binom{n}{j} { {\mathfrak{D}}_{n-j}^{(k)}(\gamma)} {\mathfrak{B}}_{j}^{[m-1, a]}(\eta;\lambda).
\end{equation}
\begin{proof}
	Using equation \eqref{e2.1}, we obtain
 \begin{equation*}
 \aligned
 \sum_{n=0}^{\infty}\frac{x^n}{n!}\, {_{\mathfrak{B}} {\mathfrak{D}}_{n, a}^{(k)}(\gamma, \eta; \lambda)} =& (1+x)^{\gamma}\,\,\frac{log(1+x)}{L_{ik}(1-e^{-x})} e^{\eta{x}}\left(\frac{x^m}{\lambda e^{x}-\sum\limits_{l=0}^{m-1}\frac{x^l}{l!}}\right)^{a}\\
 =& \left(\sum_{n=0}^{\infty}\frac{x^n}{n!}\,  {\mathfrak{D}}_{n}^{(k)}(\gamma)\right)\left(\sum_{j=0}^{\infty}\frac{x^j}{j!}\,  {\mathfrak{B}}_{j}^{[m-1, a]}(\eta; \lambda)\right)\\
  =& \sum_{n=0}^{\infty}\,\, \sum_{j=0}^{\infty}\frac{x^{n+j}}{n! j!}{\mathfrak{D}}_{n}^{(k)}(\gamma){\mathfrak{B}}_{j}^{[m-1, a]}(\eta; \lambda).
 \endaligned
 \end{equation*}
 Replceing n by n-j and and using series arrangement technique, we obtain
 
 \begin{equation*}
 \sum_{n=0}^{\infty}\frac{x^n}{n!}\, {_{\mathfrak{B}} {\mathfrak{D}}_{n, a}^{[k, m-1]}(\gamma, \eta; \lambda)} =\sum_{n=0}^{\infty}\frac{x^{n}}{n!}\left( \sum\limits_{j=0}^{n} \binom{n}{j} {\mathfrak{D}}_{n-j}^{(k)}(\gamma){\mathfrak{B}}_{j}^{[m-1, a]}(\eta; \lambda)\right).
 \end{equation*}
 
By equating both side with same power of $x^{n}$ we get desired result.
\end{proof}
\end{theorem}

\begin{theorem}
Let $n, m \in \mathbb{N}_0$, $\gamma, \eta \in \mathbb{R}$ then GABPDP satisfy following relation:
\begin{equation}\label{e2.3}
{_{\mathfrak{B}} {\mathfrak{D}}_{n, a}^{[k, m-1]}(\gamma, \eta; \lambda)} = \frac{{_{\mathfrak{B}} {\mathfrak{D}}_{n+1, a}^{[k, m-1]}(\gamma+1, \eta; \lambda)}-{_{\mathfrak{B}} {\mathfrak{D}}_{n+1, a}^{[k, m-1]}(\gamma, \eta; \lambda)}}{n+1}.
\end{equation}
\begin{proof}
Using definition \eqref{e2.1}, we have
$$\sum_{n=0}^{\infty}\frac{x^n}{n!}\, {_{\mathfrak{B}} {\mathfrak{D}}_{n, a}^{(k)}(\gamma+1, \eta; \lambda)} - \sum_{n=0}^{\infty}\frac{x^n}{n!}\, {_{\mathfrak{B}} {\mathfrak{D}}_{n, a}^{(k)}(\gamma, \eta; \lambda)}$$

$=(1+x)^{\gamma+1}\,\,\frac{log(1+x)}{L_{ik}(1-e^{-x})} e^{\eta{x}}\left(\frac{x^m}{\lambda e^{x}-\sum\limits_{l=0}^{m-1}\frac{x^l}{l!}}\right)^{a}-(1+x)^{\gamma}\,\,\frac{log(1+x)}{L_{ik}(1-e^{-x})} e^{\eta{x}}\left(\frac{x^m}{\lambda e^{x}-\sum\limits_{l=0}^{m-1}\frac{x^l}{l!}}\right)^{a}$

$=(1+x)^{\gamma}\,\,\frac{log(1+x)}{L_{ik}(1-e^{-x})} e^{\eta{x}}\left(\frac{x^m}{\lambda e^{x}-\sum\limits_{l=0}^{m-1}\frac{x^l}{l!}}\right)^{a}(1+x-1)$

\vspace{0.25cm}
$=\sum_{n=0}^{\infty}\frac{x^{n+1}}{n!}\, {_{\mathfrak{B}} {\mathfrak{D}}_{n, a}^{[k, m-1]}(\gamma, \eta; \lambda)}$

\vspace{0.25cm}
Replacing n by n+1 left hand side and equating both side the same power of $x^{n+1}$, we get desired result \eqref{e2.3}.
\end{proof}
\end{theorem}

\begin{theorem}
Let $\gamma, \eta, \omega \in \mathbb{R}$, $n, j \in \mathbb{N}_0$ and $ m \in \mathbb{N}$ then generalized Apostol-Bernoulli poly-Daehee polynomial satisfy following relation:

\begin{equation}\label{e2.4}
{_{\mathfrak{B}} {\mathfrak{D}}_{n, a}^{[k, m-1]}(\gamma+\omega, \eta; \lambda)}= \sum\limits_{j=0}^{n}\binom{n}{j} \,\,{_{\mathfrak{B}} {\mathfrak{D}}_{n-j, a}^{[k, m-1]}(\gamma, \eta; \lambda)}\,\,(\omega)_{j}
\end{equation}
\begin{proof}
	we know that
	\begin{equation*}
	(\omega)_{j}= \omega(\omega-1)(\omega-2)...(\omega-j+1)
	\end{equation*}
	By using above relation and replace $\gamma$ by $\gamma+\omega$ in equtaion \eqref{e2.1}, we get 
	\begin{equation*}
	\aligned
	\sum_{n=0}^{\infty}\frac{x^n}{n!}\, {_{\mathfrak{B}} {\mathfrak{D}}_{n, a}^{[k, m-1]}(\gamma+\omega, \eta; \lambda)} = & (1+x)^{\gamma+\omega}\,\,\frac{log(1+x)}{L_{ik}(1-e^{-x})} e^{\eta{x}}\left(\frac{x^m}{\lambda e^{x}-\sum\limits_{l=0}^{m-1}\frac{x^l}{l!}}\right)^{a}\\
	=&(1+x)^{\gamma}\,\,\frac{log(1+x)}{L_{ik}(1-e^{-x})} e^{\eta{x}}\left(\frac{x^m}{\lambda e^{x}-\sum\limits_{l=0}^{m-1}\frac{x^l}{l!}}\right)^{a}(1+x)^{\omega}\\
	=& \left(\sum_{n=0}^{\infty}\frac{x^n}{n!}\, {_{B} {\mathfrak{D}}_{n, a}^{[k, m-1]}(\gamma, \eta; \lambda)}\right)\left(\sum_{j=0}^{\infty}(\omega)_{j}\frac{x^{j}}{j!}\right)\\
	 =& \sum_{n=0}^{\infty}\,\, \sum_{j=0}^{\infty}\frac{x^{n+j}}{n! j!}{\mathfrak{D}}_{n, a}^{[k, m-1]}(\gamma, \eta; \lambda)(\omega)_{j}.
	\endaligned
	\end{equation*}
	
\vspace{0.25cm}
Replacing n by n-j and comparing both side the coefficient of $x^{n}$ of above equation, we get the desired result \eqref{e2.4}.
\end{proof}
\end{theorem}

\begin{theorem}
Let $\gamma, \eta \in \mathbb{R}$, $n, j \in \mathbb{N}_0$ and $m \in \mathbb{N}$ then GABPDP satisfy following relation:
\begin{equation}\label{e2.5}
\sum\limits_{j=0}^{n}\binom{n}{j}\,{\mathfrak{B}}_{j}^{k} {_{\mathfrak{B}} {\mathfrak{D}}_{n-j, a}^{[k, m-1]}(\gamma, \eta; \lambda)}=\sum\limits_{j=0}^{n}\binom{n}{j}\,{\mathfrak{B}}_{j}\, {_{\mathfrak{B}} {\mathfrak{D}}_{n-j, a}^{[ m-1]}(\gamma, \eta; \lambda)}.
\end{equation}
\begin{proof}
First we solve for l.h.s. by using \eqref{e1.8} and \eqref{e2.1}

$$(1+x)^{\gamma}\,\,\frac{log(1+x)}{(e^{x}-1)} e^{\eta{x}}\left(\frac{x^m}{\lambda e^{x}-\sum\limits_{l=0}^{m-1}\frac{x^l}{l!}}\right)^{a}$$
\begin{equation*}
\aligned
=&\frac{L_{ik}(1-e^{-x})}{(e^{x}-1)}(1+x)^{\gamma}\,\,\frac{log(1+x)}{L_{ik}(1-e^{-x})} e^{\eta{x}}\left(\frac{x^m}{\lambda e^{x}-\sum\limits_{l=0}^{m-1}\frac{x^l}{l!}}\right)^{a}\\
=& \left(\sum\limits_{j=0}^{\infty}\, \frac{x^j}{j!}\,{\mathfrak{B}}_{j}^{k}\right)\left(\sum_{n=0}^{\infty}\frac{x^n}{n!}\, {_{\mathfrak{B}} {\mathfrak{D}}_{n, a}^{[k, m-1]}(\gamma, \eta; \lambda)}\right)\\
=& \sum_{n=0}^{\infty}\,\, \sum_{j=0}^{\infty}\frac{x^{n+j}}{n! j!}\,{\mathfrak{B}}_{j}^{k}\,{\mathfrak{D}}_{n, a}^{[k, m-1]}(\gamma, \eta; \lambda).
\endaligned
\end{equation*}
By replacing n to n-j and arranging the terms in the the above relation, we obtain
\begin{equation}\label{e2.6}
(1+x)^{\gamma}\,\,\frac{log(1+x)}{(e^{x}-1)} e^{\eta{x}}\left(\frac{x^m}{\lambda e^{x}-\sum\limits_{l=0}^{m-1}\frac{x^l}{l!}}\right)^{a}= \sum_{n=0}^{\infty}\frac{x^{n}}{n!}\left( \sum\limits_{j=0}^{n}\binom{n}{j}\,{\mathfrak{B}}_{j}^{k}\,{\mathfrak{D}}_{n-j, a}^{[k, m-1]}(\gamma, \eta; \lambda)\right).
\end{equation}
Again, solve for r.h.s. using \eqref{e1.2} and \eqref{e1.6}, we get
$$(1+x)^{\gamma}\,\,\frac{log(1+x)}{(e^{x}-1)} e^{\eta{x}}\left(\frac{x^m}{\lambda e^{x}-\sum\limits_{l=0}^{m-1}\frac{x^l}{l!}}\right)^{a}$$
\begin{equation*}
\aligned
=&\frac{x}{(e^{x}-1)}(1+x)^{\gamma}\,\,\frac{log(1+x)}{x} e^{\eta{x}}\left(\frac{x^m}{\lambda e^{x}-\sum\limits_{l=0}^{m-1}\frac{x^l}{l!}}\right)^{a}\\
=& \left( \sum_{j=0}^{\infty}\frac{x^j}{n!}{\mathfrak{B}}_{j}\right)\left(\sum\limits_{n=0}^{\infty}{\mathfrak{D}}_{n-j, a}^{[ m-1]}(\gamma, \eta; \lambda)\right)
\endaligned
\end{equation*}

By arranging the terms of above equation and replacing n to n-j, we obtain
\begin{equation}\label{e2.7}
(1+x)^{\gamma}\,\,\frac{log(1+x)}{(e^{x}-1)} e^{\eta{x}}\left(\frac{x^m}{\lambda e^{x}-\sum\limits_{l=0}^{m-1}\frac{x^l}{l!}}\right)^{a}= \sum_{n=0}^{\infty}\frac{x^n}{n!}\left(\sum\limits_{j=0}^{n}\binom{n}{j}\,{\mathfrak{B}}_{j}\,{\mathfrak{D}}_{n-j, a}^{[ m-1]}(\gamma, \eta; \lambda)\right)
\end{equation}
Therefore, by comparing the equations \eqref{e2.6} and \eqref{e2.7}, we get the desired result \eqref{e2.5}.
\end{proof}
\end{theorem}

\begin{theorem}
	Let $\gamma, \eta, \lambda \in \mathbb{R}$, $n, j \in \mathbb{N}_0$ and $m \in \mathbb{N}$ the following relation hold true:
\begin{equation}\label{e2.8}
{\mathfrak{B}}_{n, a}^{[m-1]}(\eta; \lambda)=\sum\limits_{j=0}^{n}\binom{n}{j}\,b_{j}^{k}(-\gamma)\, {_{\mathfrak{B}} {\mathfrak{D}}_{n-j, a}^{[k, m-1]}(\gamma, \eta; \lambda)}.
\end{equation}
\begin{proof}
 By using \eqref{e1.22} and \eqref{e2.1}, 
 \begin{equation*}
 \aligned
\sum_{n=0}^{\infty}\frac{x^n}{n!}{\mathfrak{B}}_{n, a}^{[m-1]}(\eta; \lambda)=& e^{\eta{x}}\left(\frac{x^m}{\lambda e^{x}-\sum\limits_{l=0}^{m-1}\frac{x^l}{l!}}\right)^{a}\\
=&(1+x)^{-\gamma}\left(\frac{L_{ik}(1-e^{-x})}{log(1+x)}\right)\sum_{n=0}^{\infty}\frac{x^n}{n!}\, {_{\mathfrak{B}} {\mathfrak{D}}_{n, a}^{[k, m-1]}(\gamma, \eta; \lambda)}\\
=& \left(\sum_{j=0}^{\infty}\frac{x^j}{j!}{\mathfrak{b}}_{j}^{(k)}(-\gamma)\right)\left(\sum_{n=0}^{\infty}\frac{x^n}{n!}\, {_{\mathfrak{B}} {\mathfrak{D}}_{n, a}^{[k, m-1]}(\gamma, \eta; \lambda)}\right)
\endaligned
\end{equation*}
Therefore, by using series arrangement technique and equating the coefficient of $x^{n}$ both side in the above relation, we get required result \eqref{e2.8}.
\end{proof}
\end{theorem}

\section{\bf{Implicit and summation formulae for GABPDP}}

In this section, we drive a some useful implicit and summation formula of generalized Apostol-Bernoulli based poly daehee polynomial (GABPDP).

\begin{theorem}
Let $\gamma, \eta, \lambda \in \mathbb{R}, n, p, q, b, c \in \mathbb{N}_0, m \in \mathbb{N}$ and $a \in \mathbb{C}$, the following implicit summation formula for generalized Apostol-Bernoulli poly-Daehee polynomials holds true:
\begin{equation}\label{e3.1}
_{\mathfrak{B}} {\mathfrak{D}}_{n, a, b+c}^{[k, m-1]}(\gamma, \eta; \lambda)= \sum\limits_{n, q=0}^{b, c}\binom{b}{n}\binom{c}{q}(\eta-\omega)^{n+q}~~_{\mathfrak{B}} {\mathfrak{D}}_{n, a; b+c-n-q}^{[k, m-1]}(\gamma, \omega; \lambda)
\end{equation}
\begin{proof}
We know following series manipulation formulae	
\begin{equation}\label{E3.2}
\sum_{N=0}^{\infty} g(N)\frac{(\gamma + \eta)^N}{N!}= \sum_{n, m = 0}^{\infty}\frac{\gamma^{n}}{n!}\frac{\eta^{m}}{m!}g(n+m)
\end{equation}
Now, replacing $x$ by $x+\mu$ in \eqref{e2.1}, we get
$$(1+(x+\mu))^{\gamma}\,\,\frac{log(1+(x+\mu))}{L_{ik}(1-e^{-(x+\mu)})}~~ e^{\eta{(x+\mu)}}\left(\frac{(x+\mu)^m}{\lambda e^{(x+\mu)}-\sum\limits_{l=0}^{m-1}\frac{(x+\mu)^l}{l!}}\right)^{a}$$
\begin{equation}\label{E3.3}
=\sum_{b, c = 0}^{\infty}\frac{x^{b}}{b!}\frac{\mu^{c}}{c!}\, {_{\mathfrak{B}} {\mathfrak{D}}_{n, a; b+c}^{[k, m-1]}(\gamma, \eta; \lambda)}
\end{equation}
Again, replacing $\eta$ by $\omega$ in \eqref{E3.3} and equate with \eqref{E3.3}, we have
\begin{equation}\label{E3.4}
e^{(\eta-\omega)(x+\mu)}\sum_{b, c =0}^{\infty}\frac{x^{b}}{b!}\frac{\mu^{c}}{c!}\, {_{\mathfrak{B}} {\mathfrak{D}}_{n, a; b+c}^{[k, m-1]}(\gamma, \omega; \lambda)}=\sum_{b, c = 0}^{\infty}\frac{x^{b}}{b!}\frac{\mu^{c}}{c!}\, {_{\mathfrak{B}} {\mathfrak{D}}_{n, a; b+c}^{[k, m-1]}(\gamma, \eta; \lambda)}
\end{equation}
\begin{equation}\label{E3.5}
\sum_{N=0}^{\infty}\frac{[(\eta-\omega)(x+\mu)]^{N}}{N!}\sum_{b, c= 0}^{\infty}\frac{x^{b}}{b!}\frac{\mu^{c}}{c!}\, {_{\mathfrak{B}} {\mathfrak{D}}_{n, a; b+c}^{[k, m-1]}(\gamma, \omega; \lambda)}=\sum_{b, c= 0}^{\infty}\frac{x^{b}}{b!}\frac{\mu^{c}}{c!}\, {_{\mathfrak{B}} {\mathfrak{D}}_{n, a; b+c}^{[k, m-1]}(\gamma, \eta; \lambda)}.
\end{equation}
Using \eqref{E3.2} in the l.h.s. of \eqref{E3.5}, we get
\begin{equation}\label{E3.6}
\sum\limits_{n, q=0}^{b, c}\frac{(\eta - \omega)^{n+q}x^{n} \mu^{q}}{n!q!}\sum_{b, c=0}^{\infty}\frac{x^{b}}{b!}\frac{\mu^{c}}{c!}\, {_{\mathfrak{B}} {\mathfrak{D}}_{n, a; b+c}^{[k, m-1]}(\gamma, \omega; \lambda)}=\sum_{b, c= 0}^{\infty}\frac{x^{b}}{b!}\frac{\mu^{c}}{c!}\, {_{\mathfrak{B}} {\mathfrak{D}}_{n, a; b+c}^{[k, m-1]}(\gamma, \eta; \lambda)}
\end{equation}
Using l.h.s. series arrangement method after that replacing $b$ by $b-n$  and $c$ by $c-q$ in \eqref{E3.6}, we have
$$\sum_{b, c=0}^{\infty}\sum\limits_{n, q=0}^{b, c}\frac{(\eta - \omega)^{n+q}}{n!q!}\, {_{\mathfrak{B}} D_{n, a; b+c-n-q}^{[k, m-1]}(\gamma, \omega; \lambda)}\frac{x^{b}}{(b-n)!}\frac{\mu^{c}}{(c-q)!}$$
\begin{equation}\label{E3.7}
=\sum_{b, c= 0}^{\infty}\frac{x^{b}}{b!}\frac{\mu^{c}}{c!}\, {_{\mathfrak{B}} {\mathfrak{D}}_{n, a; b+c}^{[k, m-1]}(\gamma, \eta; \lambda)}
\end{equation}
Therefore, by equating the coefficient of same power of $x^{b}$ and $x^{c}$ in \eqref{E3.7}, we get the desired result.
\end{proof}
\end{theorem}

\begin{theorem}
For $\gamma, \eta,\omega, \lambda \in \mathbb{R}$, $n \in \mathbb{N}_0$, $m \in \mathbb{N}$ and $a, b \in \mathbb{C}$, then GABPDP satisfy the following relation:
\begin{equation}\label{e3.2}
 {_{\mathfrak{B}} {\mathfrak{D}}_{n, a+b}^{[k, m-1]}(\gamma, \eta+\omega; \lambda)}= \sum\limits_{j=0}^{n}\binom{n}{j}\,\,_{\mathfrak{B}} {\mathfrak{D}}_{n-j, a}^{[k, m-1]}(\gamma, \eta; \lambda)\,{\mathfrak{B}}_{j, b}^{[m-1]}(\omega, \lambda).
\end{equation}

\begin{proof}
Replacing $\eta$ by $\eta+\omega$ and $a$ by $a+b$ in equation \eqref{e2.1}, we obtain
\begin{equation*}
\sum_{n= 0}^{\infty}\frac{x^n}{n!}\, {_{\mathfrak{B}} {\mathfrak{D}}_{n, a+b}^{[k, m-1]}(\gamma, \eta+\omega; \lambda)}=(1+x)^{\gamma}\,\,\frac{log(1+x)}{L_{ik}(1-e^{-x})}\,\, e^{(\eta+\omega){x}}\left(\frac{x^m}{\lambda e^{x}-\sum\limits_{l=0}^{m-1}\frac{x^l}{l!}}\right)^{a+b}.
\end{equation*}
\begin{equation*}
\aligned
=&\left((1+x)^{\gamma}\,\,\frac{log(1+x)}{L_{ik}(1-e^{-x})}\,\, e^{\eta{x}}\left(\frac{x^m}{\lambda e^{x}-\sum\limits_{l=0}^{m-1}\frac{x^l}{l!}}\right)^{a}\right)\left(e^{\omega{x}}\left(\frac{x^m}{\lambda e^{x}-\sum\limits_{l=0}^{m-1}\frac{x^l}{l!}}\right)^{b}\right)\\
=& \left(\sum_{n= 0}^{\infty}\frac{x^n}{n!}\, {_{\mathfrak{B}} {\mathfrak{D}}_{n, a}^{[k, m-1]}(\gamma, \eta; \lambda)}\right)\left(\sum_{j= 0}^{\infty}\frac{x^j}{j!}{\mathfrak{B}}_{j, b}^{[m-1]}(\omega; \lambda)\right).
\endaligned
\end{equation*}
First, we apply a series arrangement method after that comparing both sides to $x^{n}$, we get required result \eqref{e3.2}
\end{proof}
\end{theorem}

\begin{theorem}
For $\gamma, \eta \in \mathbb{R}$, $n \in \mathbb{N}_0, m \in \mathbb{N}$ and $a \in \mathbb{C}$, the generalized Apostol-Bernoulli poly-Daehee Polynomials satisfy following relation:
\begin{equation}\label{e3.3}
{_{\mathfrak{B}} {\mathfrak{D}}_{n, a}^{[k, m-1]}(\gamma, \eta; \lambda)}=\sum\limits_{j=0}^{n}\binom{n}{j}\,\, {\mathfrak{D}}_{n-j, a}^{[k, m-1]}(\gamma; \eta-\omega)\,{\mathfrak{B}}_{j, b}^{[m-1]}(\omega, \lambda)
\end{equation}
\begin{proof}
By using \eqref{e2.1}, we can write
$$(1+x)^{\gamma}\,\,\frac{log(1+x)}{L_{ik}(1-e^{-x})}\,\, \left(\frac{x^m}{\lambda e^{x}-\sum\limits_{l=0}^{m-1}\frac{x^l}{l!}}\right)^{a}~~e^{(\eta-\omega){x}+\omega{x}}$$
\begin{equation*}
\aligned
=& \left((1+x)^{\gamma}\,\,\frac{log(1+x)}{L_{ik}(1-e^{-x})}e^{(\eta-\omega){x}}\right)\left(e^{\omega{x}}\left(\frac{x^m}{\lambda e^{x}-\sum\limits_{l=0}^{m-1}\frac{x^l}{l!}}\right)^{a}\right)\\
=& \left(\sum_{n=0}^{\infty}{\mathfrak{D}}_{n, a}^{[k, m-1]}(\gamma; \eta-\omega)\right)\left(\sum_{j=0}^{\infty}{\mathfrak{B}}_{j, b}^{[m-1]}(\omega; \lambda)\right)
\endaligned
\end{equation*}
By using series arrangement technique and comparing both sides $x^{n}$, we get desired result \eqref{e3.3}.
\end{proof}
\end{theorem}

\begin{theorem}
Let $\gamma, \eta, \lambda \in \mathbb{R}, n, j \in \mathbb{N}_0, m \in \mathbb{N}$ and $a \in \mathbb{C}$, then GABPDP satisfy following relation:
\begin{equation}\label{e3.4}
{_{\mathfrak{B}} {\mathfrak{D}}_{n, a}^{[k, m-1]}(\gamma, \eta+1; \lambda)}= \sum\limits_{j=0}^{n}\binom{n}{j}\,\,_{\mathfrak{B}} {\mathfrak{D}}_{n-j, a}^{[k, m-1]}(\gamma, \eta; \lambda)
\end{equation}
\begin{proof}
Using \eqref{e2.1} and replacing $\eta$ by $\eta+1$, we obtain 
\begin{equation*}
\aligned
\sum_{n=0}^{\infty}\frac{x^n}{n!}\, {_{\mathfrak{B}} {\mathfrak{D}}_{n, a}^{[k, m-1]}(\gamma, \eta+1; \lambda)}=&(1+x)^{\gamma}\,\,\frac{log(1+x)}{L_{ik}(1-e^{-x})} \left(\frac{x^m}{\lambda e^{x}-\sum\limits_{l=0}^{m-1}\frac{x^l}{l!}}\right)^{a}~~e^{(\eta+1){x}}\\
=&\left((1+x)^{\gamma}\,\,\frac{log(1+x)}{L_{ik}(1-e^{-x})} \left(\frac{x^m}{\lambda e^{x}-\sum\limits_{l=0}^{m-1}\frac{x^l}{l!}}\right)^{a}~e^{\eta{x}}\right)e^{x}\\
=&\left(\sum_{n=0}^{\infty}\frac{x^n}{n!}\, {_{\mathfrak{B}} {\mathfrak{D}}_{n, a}^{[k, m-1]}(\gamma, \eta; \lambda)}\right)\left(\sum\limits_{j=0}^{\infty}\frac{x^j}{j!}\right).
\endaligned
\end{equation*}
By using series arrangement technique and comparing both side $x^{n}$, we get required result \eqref{e3.4}.
\end{proof}
\end{theorem}

\begingroup
\begin{table}
	\caption {\label{tab:table1} {\bf Members similar to the polynomials ${_{\mathfrak{B}} {\mathfrak{D}}_{n, a}^{[k, m-1]}(\gamma, \eta; \lambda)}$}.} 
	\scriptsize
	\begin{tabular}{llllll}
		
		&  {\bf S. No.}   & {\bf Name of polynomial} & 
		&  $\bf {\mathcal{A}(t)}$ &                                                            
		{\bf Generating function}                                                  \\    \hline
		&I.   & Apostol-Bernoulli based \\ & &poly-Daehee polynomials                                           &
		& $\left(\frac{x}{\lambda e^{x}-1}\right)^{a}\,\frac{\log(1+x)}{L_{ik}(1-e^{-x})}(1+x)^{\gamma}$ & 
		$\sum\limits_{n=0}^{\infty}{}_{\mathfrak{B}} {\mathfrak{D}}^{(k)}_{n,\alpha}(\gamma,\eta;\lambda)\frac{x^{n}}{n!}$\\ & & & & & $=\left(\frac{2}{\lambda e^{x}+1}\right)^{a}\,\frac{\log(1+x)}{L_{ik}(1-e^{-x})}(1+x)^{\gamma}\,e^{\eta{x}}$         \\
		
		&II.   & Apostol-Euler based \\ & &poly-Daehee polynomials                                         &
		& $\left(\frac{2}{\lambda e^{x}+1}\right)^{a}\,\frac{\log(1+x)}{L_{ik}(1-e^{-x})}(1+x)^{\gamma}$ & 
		$\sum\limits_{n=0}^{\infty}{}_\mathcal{E} {\mathfrak{D}}^{(k)}_{n,\alpha}(\gamma,\eta;\lambda)\frac{x^{n}}{n!}$\\ & & & & & $=\left(\frac{2}{\lambda e^{x}+1}\right)^{a}\,\frac{\log(1+x)}{L_{ik}(1-e^{-x})}(1+x)^{\gamma}\,e^{\eta{x}}$         \\ 
		
		&III.   & Apostol-Genocchi based \\ & &poly-Daehee polynomials                                     &
		& $\left(\frac{2x}{\lambda e^{x}-1}\right)^{a}\frac{\log(1+x)}{L_{ik}(1-e^{-x})}(1+x)^{\eta}$& 
		$\sum\limits_{n=0}^{\infty}{}_\mathcal{G} {\mathfrak{D}}^{(k)}_{n,\alpha}(\gamma,\eta;\lambda)\frac{x^{n}}{n!}$\\ & & & & &  $=\left(\frac{2x}{\lambda e^{x}-1}\right)^{a}\,\frac{\log(1+x)}{L_{ik}(1-e^{-x})}(1+x)^{\gamma}\,e^{\eta{x}}$         \\ 
		\\         
		&IV.   & Apostol-Bernoulli based \\ & &Daehee polynomials                                   &
		& $\left(\frac{x}{\lambda e^{x}-1}\right)^{a}\frac{\log(1+x)}{x}(1+x)^{\gamma}$
		&  $\sum\limits_{n=0}^{\infty}{}_{\mathfrak{B}} {\mathfrak{D}}_{n, a}(\gamma,\eta;\lambda)\frac{x^{n}}{n!}$\\ & & & & & $=\left(\frac{x}{\lambda e^{x}-1}\right)^{\alpha}\,\frac{\log(1+x)}{x}(1+x)^{\gamma}\,e^{\eta{x}}$         \\ 
		
		\\
		&V.   & Apostol-Euler based \\ & &Daehee polynomials                                    &
		& $\left(\frac{2}{\lambda e^{x}+1}\right)^{a}\frac{\log(1+x)}{x}(1+x)^{\gamma}$ & $\sum\limits_{n=0}^{\infty}{}_\mathcal{E}{\mathfrak{D}}_{n, a}(\gamma,\eta;\lambda)\frac{x^{n}}{n!}$\\ & & & & & $=\left(\frac{2}{\lambda e^{x}+1}\right)^{a}\,\frac{\log(1+x)}{x}(1+x)^{\gamma}\,e^{\eta{x}}$         \\ 
		
		\\
		&VI.   & Apostol-Genocchi based \\ & &Daehee polynomials                                    &
		& $\left(\frac{2x}{\lambda e^{x}+1}\right)^{a}\frac{\log(1+x)}{x}(1+x)^{\gamma}$ & $=\sum\limits_{n=0}^{\infty}{}_\mathcal{G}{\mathfrak{D}}_{n, a}(\gamma,\eta;\lambda)\frac{x^{n}}{n!}$\\ & & & & & $=\left(\frac{2x}{\lambda e^{x}+1}\right)^{a}\,\frac{\log(1+x)}{x}(1+x)^{\gamma}\,e^{\eta{x}}$         \\ 
		\hline          \\    
		\\
	\end{tabular}
	
\end{table}
\endgroup 
\newpage

\section{\bf{Conclusion}}

In this paper, we defined the relation between generalized Apostol-Bernoulli polynomial and poly Daehee polynomial known as generalized Apostol-Bernoulli poly-Daehee polynomial (GABPDP) with generating function \eqref{e2.1}. Which is very useful to us, because polynomial plays a very important role to use as a solution of different kind of differential equations, also polynomial use to represent a different kind of characteristic linear dynamic system. We also discus some useful identities and their implicit summation formulae of generalized Apostol-Bernoulli poly-Daehee polynomials (GABPDP).


\begin{thebibliography}{99}

\bibitem{carlitz} L. Carlitz, {\it A note on Bernoulli and Euler polynomials of the second kind}, Scripta Math 25 (1961): 323-330.

\bibitem{choi-nu-usman1} J. Choi, N.U. Khan, T. Usman,  and M. Aman,  {\it Certain unified polynomials}, Integral Transforms and Special Functions, 30(1) (2019): 28-40 .

\bibitem{Dattoli} G. Cesarano, S. Dattoli and C. Lorenzutta {\it Finite sums and generalized forms of Bernoulli polynomials}, Rendiconti di Matematica 19 (1999): 385-391.

\bibitem{jolany-mohsen} J. Hassan, M. Aliabadi, R.B. Corcino and M.R. Darafsheh {\it A note on multi poly-Euler numbers and Bernoulli polynomials}, arXiv preprint arXiv:1401.2645 (2014).

\bibitem{khan-usman-aman} N.U. Khan, T. Usman, M. Aman  {\it Generating function for Legendre-based poly-Bernoulli numbers and polnomials}, Honam Math. J.  39 (2017): 217-231.

\bibitem{khan-usman1} N. U. Khan, and T. Usman {\it A new class of Laguerre-based generalized Apostol polynomials} Fasciculi Mathematici 57, no. 1 (2016): 67-89.

\bibitem{khan-usman2}N. U. Khan, and  T. Usman {\it A new class of Laguerre poly-Bernoulli numbers and polynomials}, Adv Stud Contemporary Math 27 (2017): 229-241.

\bibitem{khan-usman3} N. U. Khan, and T. Usman {\it A new class of Laguerre-based poly-Euler and multi poly-Euler polynomials}, J Anal Num Theor 4 (2016): 113-120.

\bibitem{khan-usman-choi2} N.U. Khan, T. Usman,  and  J. Choi, {\it A new class of generalized polynomials}, Turkish Journal of Mathematics 42, no. 3 (2018): 1366-1379.

\bibitem{wasim}
W.A. Khan  {\it A new class of Hermite poly-Genocchi polynomials}, J. Anal. Number Theory,  (4) (2016): 1-8.

\bibitem{Waseem} W. A.Khan {\it A note on Hermite-based poly-Euler and multi poly-Euler polynomials}, Palestine J. Math 5, no. 1 (2016): 17-26.

\bibitem{subuhi-tabinda}S. Khan, and T. Nahid {\it Finding non-linear differential equations and certain identities for the Bernoulli–Euler and Bernoulli–Genocchi numbers}, SN Applied Sciences 1, no. 3 (2019): 217.

\bibitem{san-taekyun} D.S. Kim and  T. Kim, {\it Some identities involving Genocchi polynomials and numbers}, ARS COMBINATORIA 121 (2015): 403-412.

\bibitem{kim-taekyun}  T. Kim {\it On the multiple q-Genocchi and Euler numbers}, Russian Journal of Mathematical Physics 15, no. 4 (2008): 481-486.

\bibitem{kim-kwon-lee} T. Kim, H.I. Kwon,  S.H. Lee and J.J. Seo  {\it A note on poly-Bernoulli numbers and polynomials of the second kind}, Advances in Difference Equations 2014, no. 1 (2014): 219.

\bibitem{dae-taekyun} D.S. Kim,  and  T. Kim {\it A study on the integral of the product of several Bernoulli polynomials}, The Rocky Mountain Journal of Mathematics 44, no. 4 (2014): 1251-1263.

\bibitem{Takao-luca} T. Komatsu and T. Luca {\it Some relationships between poly-Cauchy numbers and poly-Bernoulli numbers}, In Annales Mathematicae et Informaticae, pp. 99-105. 2013.

\bibitem{lim-kwon}
 D.S. Lim, and J. Kwon {\it A note on poly-Daehee numbers and polynomials}, {Proc. Jangjeon Math. Soc.}, \textbf{19} (2), 2016, 219--224.

\bibitem{Luo-ming-sri} Q.M. Luo, and H. M. Srivastava {\it Some generalizations of the Apostol-Bernoulli and Apostol Euler polynomials}, Journal of Mathematical Analysis and Applications 308, no. 1 (2005): 290-302.

\bibitem{luo-ming} Q.M. Luo  {\it Apostol Euler polynomials of higher order and Gaussian hypergeometric functions}, Taiwanese Journal of Mathematics 10, no. 4 (2006): 917-925.

\bibitem{luo-ming2} Q.M. Luo {\it Extensions of the Genocchi polynomials and their Fourier expansions and integral representations} Osaka Journal of Mathematics 48, no. 2 (2011): 291-309.

\bibitem{Rain} E.D. Rainville {\it Special Functions}, The MacMillan Comp., New York (1960)

\bibitem{Tremblay-Gaboury}T. Richard, S. Gaboury and B.J. Fugere {\it A new class of generalized Apostol-Bernoulli polynomials and some analogues of the Srivastava–Pintér addition theorem}, Applied Mathematics Letters 24, no. 11 (2011): 1888-1893.



\end{thebibliography}
\end{document}